\documentclass[12pt,twoside]{article}

\def\Z{{\rlap{$\kern2pt{\rm Z}$}{\rm Z}\,}}
\def\bld#1#2{{\buildrel{#1}\over{#2}}}
\def\st#1#2{{\mathrel{\mathop{#2}\limits_{#1}}{}\!}}
\def\stb#1#2#3{{\st{{#1}}{\bld{{#2}}{#3}}{}\!}}
\def\xmare#1#2{\stb{#1}{#2}{\mbox{\Huge$\times$}}}
\def\bldl#1{{\buildrel{#1}\over{{\longrightarrow}}}}

\def\Im{{\rm Im}}
\def\Fix{{\rm Fix}}
\def\slat{{${\cal E}$--lattice}}
\def\slats{{${\cal E}$--lattices}}
\def\Ker{{\rm Ker}}
\def\n{\noindent }
\def\Fixe{{{\rm Fix}\,{\varepsilon}}}

\def\w{\wedge }

\def\we{\wedge_{\varepsilon} }
\def\ve{\vee_{\varepsilon} }

\title{On isomorphisms of canonical ${\cal E}$-lattices}
\author{Marius T\u arn\u auceanu\\
{\small Faculty of  Mathematics, "Al.I. Cuza" University,
Ia\c si, Romania}\\
{\small e-mail: tarnauc@uaic.ro}}
\date{October 1, 2007}

\begin{document}
\maketitle

\begin{abstract} The aim of the present paper is to study isomorphisms of canonical ${\cal E}$-lattices. Some interesting results are obtained in the particular case of isomorphisms between two subgroup ${\cal E}$-lattices.
\end{abstract}
\bigskip
\bigskip

\section{Preliminaries}

The starting point for our discussion is given by the paper [10],
where there is introduced the category of \slats\ and there are
made some basic constructions in this category. Given a nonvoid
set $L$ and a map ${\varepsilon}:L\to L$, we denote by
$\Ker\,{\varepsilon}$ the kernel of ${\varepsilon}$ (i.e.
$\Ker\,{\varepsilon}=\{(a,b)\in L\times
L\mid{\varepsilon}(a)={\varepsilon}(b)\})$, by
$\Im\,{\varepsilon}$ the image of ${\varepsilon}$ (i.e.
$\Im\,{\varepsilon}=\{{\varepsilon}(a)\mid a\in L\})$ and by
$\Fix\,{\varepsilon}$ the set consisting of all fixed points of
${\varepsilon}$ (i.e. $\Fix\,\varepsilon=\{a\in L\mid
{\varepsilon}(a)=a\}).$ We say that $L$ is an {\it \slat}
(relative to ${\varepsilon})$ if there exist two binary operations
$\w_{\varepsilon},\vee_{\varepsilon}$ on $L$ which satisfy the
following properties:
\begin{itemize}
\item[a)] $a\w_{\varepsilon}(b\w_{\varepsilon} c)=(a\w_{\varepsilon} b) \w_{\varepsilon} c,\
a\vee_{\varepsilon}(b\vee_{\varepsilon} c)=(a\vee_{\varepsilon} b) \vee_{\varepsilon} c,$ for all $a,b,c\in L;$
\item[b)]
$a\w_{\varepsilon} b =b\w_{\varepsilon} a,\ a\vee_{\varepsilon} b =b\vee_{\varepsilon} a,$ for all $a,b\in L;$
\item[c)] $a\w_{\varepsilon} a =a\vee_{\varepsilon} a={\varepsilon}(a),$ for any $a\in L;$
\item[d)] $a\w_{\varepsilon}(a\vee_{\varepsilon} b)=a\vee_{\varepsilon}(a\w_{\varepsilon} b)={\varepsilon}(a),$ for all $a,b\in L.$
\end{itemize}

Clearly, in an \slat\ $L$ (relative to ${\varepsilon}$) the map
${\varepsilon}$ is idempotent and
$\Im\,{\varepsilon}=\Fix\,{\varepsilon}$. Moreover, the set
$\Fix\,{\varepsilon}$ is closed under the binary operations
$\w_{\varepsilon},\vee_{\varepsilon}$ and, denoting by
$\w_{\varepsilon}^\circ,\vee_{\varepsilon}^\circ$ the restrictions
of $\w_{\varepsilon},\vee_{\varepsilon}$ to $\Fix\,{\varepsilon}$,
we have that
$(\Fix\,{\varepsilon},\w_{\varepsilon}^\circ,\vee_{\varepsilon}^\circ)$
is a lattice. The connection between the \slat\ concept and the
lattice concept is very powerful. So, if
$(L,\w_{\varepsilon},\vee_{\varepsilon})$ is an \slat\ and $\sim$
is an equivalence relation on $L$ such that
$\sim\,\subseteq\Ker\,{\varepsilon}$, then the factor set $L/\sim$
is a lattice isomorphic to the lattice $\Fix\,{\varepsilon}$.
Conversely,  if $L$ is a nonvoid set and $\sim$ is an equivalence
relation on $L$ having the property\ that the factor set $L/\sim$
is a lattice, then the set $L$ can be endowed with an \slat\
structure (relative to a map ${\varepsilon}:L\longrightarrow  L)$
such that $\sim\,\subseteq \Ker\,{\varepsilon}$ and $L/\sim\,\cong
\Fix\,{\varepsilon}$.

We say that an {\it\slat} $(L,\we,\ve)$  is a {\it canonical
\slat} if $a\we b,a\ve b\in\Fixe,$ for all $a,b\in L.$ Three
fundamental types of canonical \slats\ have been identified in
[10]. One of the most important example of a canonical \slat\ is
constituted by the lattice $L(G)$ of all subgroups of a group $G$
(called the {\it subgroup \slat} of $G)$. Here the map
${\varepsilon}$ is defined by ${\varepsilon}(H)=H_G$ (where $H_G$
is the core of $H$ in $G)$, for any $H\in L(G)$ and the binary
operations $\we,\ve$ are defined by $H_1\we
H_2={\varepsilon}(H_1)\cap{\varepsilon}(H_2)$, $H_1\ve
H_2={\varepsilon}(H_1){\varepsilon}(H_2),$ for all $H_1,H_2\in
L(G)$. Mention that in this situation the lattice $\Fixe$ is just
the normal subgroup lattice $N(G)$ of $G$.

Let  $(L_1,\wedge_{{\varepsilon}_1},\vee_{{\varepsilon}_1})$ and
$(L_2,\wedge_{{\varepsilon}_2},\vee_{{\varepsilon}_2})$ be two
\slats. A map $f:L_1\to L_2$ is called an {\it\slat\ homomorphism}
if:
\begin{itemize}
\item[a)] $f\circ{\varepsilon}_1={\varepsilon}_2\circ f;$
\item[b)] for all $a,b\in L_1$, we have:
\begin{itemize}
\item[i)] $f(a\w_{{\varepsilon}_1}b)=f(a)\w_{{\varepsilon}_2}f(b);$
\item[i)] $f(a\vee_{{\varepsilon}_1}b)=f(a)\vee_{{\varepsilon}_2}f(b).$
\end{itemize}\end{itemize}
Moreover, if the map $f$ is one-to-one and onto, then we say that
it is an {\it\slat\  isomorphism}. An \slat\ isomorphism of an
\slat\ into itself is called an {\it\slat\ automorphism}. For an
\slat\ $L$, we shall denote by ${\rm Aut}_{\cal E}(L)$ the group
consisting of all \slat\ automorphisms of $L$.

\section{Main results}

\subsection{Isomorphisms of canonical \slats}

In this section we present some general results concerned to isomorphisms between canonical \slats.

 Let $(L_1,\w_{{\varepsilon}_1},\vee_{{\varepsilon}_1})$ and
$(L_2,\w_{{\varepsilon}_2},\vee_{{\varepsilon}_2})$ be two canonical \slats. For every element $a\in L_i$, we denote by $[a]_i$ the equivalence class of $a$ modulo $\Ker\,{\varepsilon}_i$ (i.e. $[a]_i=\{b\in L_i\mid {\varepsilon}_i(b)={\varepsilon}_i(a)\})$, $i=1,2.$ First of all, we give a characterization of \slat\ isomorphisms from $L_1$ to $L_2$.

\bigskip\n{\bf Proposition 1.} {\it Let $(L_1,\wedge_{{\varepsilon}_1},\vee_{{\varepsilon}_1})$ and $(L_2,\wedge_{{\varepsilon}_2},\vee_{{\varepsilon}_2})$ be
two canonical \slats\ and  $f:L_1\to L_2$ be a map. Then the
following two conditions are equivalent:
\begin{itemize}
\item[{\rm a)}] $f$ is an \slat\ isomorphism.
\item[{\rm b)}] {\rm i)} The restriction $f_0$ of $f$ to the set $\Fix\,{\varepsilon}_1$ is a lattice isomorphism from\\\ \hspace{4mm} $\Fix\,{\varepsilon}_1$ to $\Fix\,{\varepsilon}_2$.\\
{\rm ii)} $f\mid[a]_1:[a]_1\longrightarrow [f(a)]_2$ is one-to-one and onto, for each $a\in\Fix\,{\varepsilon}_1.$
\end{itemize}}

\n{\bf Proof.} a) $\Longrightarrow $ b) Let $a\in\Fix\,{\varepsilon}_1.$ Then ${\varepsilon}_1(a)=a$ and so $f_0(a)=f(a)=f({\varepsilon}_1(a))=(f\circ{\varepsilon}_1)(a)=({\varepsilon}_2\circ f)(a)={\varepsilon}_2(f(a))\in\Fix\,{\varepsilon}_2.$ Thus $\Im\,f_0\subseteq\Fix\,{\varepsilon}_2.$ For all $a,b\in L_1$, we have $f_0(a\w^\circ_{{\varepsilon}_1}b)=f(a\w_{{\varepsilon}_1}b)=f(a)\w_{{\varepsilon}_2}f(b)=f_0(a)\w_{{\varepsilon}_2}f_0(b)=f_0(a)\w^\circ_{{\varepsilon}_2}f_0(b)$ and $
 f_0(a\vee^\circ_{{\varepsilon}_1}b)=f(a\vee_{{\varepsilon}_1}b)=f(a)\vee_{{\varepsilon}_2}f(b)=f_0(a)\vee_{{\varepsilon}_2}f_0(b)=f_0(a)\vee^\circ_{{\varepsilon}_2}f_0(b)$, which show that $f_0$ is a lattice homomorphism. Since $f$ is one-to-one and onto, it is clear that $f_0$ has the same properties\ and hence i) holds.

 Now, let $b\in[a]_1.$ Then ${\varepsilon}_1(b)={\varepsilon}_1(a)=a,$ which implies that ${\varepsilon}_2(f(b))=({\varepsilon}_2\circ f)(b)=(f\circ{\varepsilon}_1)(b)=f({\varepsilon}_1(b))=f(a).$ Therefore $\Im(f\mid[a]_1)\subseteq[f(a)]_2.$ As $f$ is one-to-one, $f\mid[a]_1$ is one-to-one, too. If $c\in[f(a)]_2,$ then ${\varepsilon}_2(c)=f(a)$ and so ${\varepsilon}_1(f^{-1}(c))=({\varepsilon}_1\circ f^{-1})(c)=(f^{-1}\circ{\varepsilon}_2)(c)=f^{-1}({\varepsilon}_2(c))=f^{-1}(f(a))=a.$ It results that $f^{-1}(c)\in[a]_1.$ Hence $f\mid[a]_1$ is onto.

 b) $\Longrightarrow $ a) For any $a\in L_1$, we have $a\in[{\varepsilon}_1(a)]_1.$ From the condition ii) of b), it obtains $f(a)\in[f({\varepsilon}_1(a))]_2,$ which implies that $({\varepsilon}_2\circ f)(a)={\varepsilon}_2(f(a))=f({\varepsilon}_1(a))=(f\circ{\varepsilon}_1)(a).$ Thus $f\circ{\varepsilon}_1={\varepsilon}_2\circ f$. Since both $L_1$ and $L_2$ are canonical \slats, we have $f(a\w_{{\varepsilon}_1}b)=f({\varepsilon}_1(a\w_{{\varepsilon}_1}b))=f({\varepsilon}_1(a)\w_{{\varepsilon}_1}{\varepsilon}_1(b))=f_0({\varepsilon}_1(a)\w^\circ_{{\varepsilon}_1}{\varepsilon}_1(b))=f_0({\varepsilon}_1(a))\w^\circ_{\varepsilon_2}f_0({\varepsilon}_1(b))=(f\circ{\varepsilon}_1)(a)\w_{{\varepsilon}_2}(f\circ{\varepsilon}_1)(b)=({\varepsilon}_2\circ f)(a)\w_{{\varepsilon}_2}({\varepsilon}_2\circ f)(b)={\varepsilon}_2(f(a))\w_{{\varepsilon}_2}{\varepsilon}_2(f(b))={\varepsilon}_2(f(a)\w_{{\varepsilon}_2}f(b))=f(a)\w_{{\varepsilon}_2}f(b)$ and, in the same manner, $f(a\vee_{{\varepsilon}_1}b)=f(a)\vee_{{\varepsilon}_2}f(b)$, for all $a,b\in L_1.$ Therefore $f$ is an \slat\ homomorphism. Clearly, the map $f$ is one-to-one and onto and hence our proof if finished.\bigskip

As a consequence of the above proposition, it obtains the next characterization of \slat\ autmorphisms of a canonical \slat.

\bigskip\n{\bf Corollary.} {\it Let $(L,\we,\ve)$ be a canonical \slat, ${\rm Aut}(\Fixe)$ be the group consisting of all automorphisms of the lattice $\Fixe$ and $f:L\longrightarrow  L$ be a map. Then $f\in{\rm Aut}_{\cal E}(L)$ if and only if $f\mid\Fix\,{\varepsilon}\in{\rm Aut}(\Fixe)$ and $f\mid[a]:[a]\longrightarrow [f(a)]$ is one-to-one and onto, for each $a\in\Fixe.$}\bigskip

In the following we shall investigate the structure of the group
${\rm Aut}_{\cal E}(L)$ as\-so\-cia\-ted to a canonical \slat\
$(L,\we,\ve)$. Suppose that \mbox{$\Fixe{=}\{a_i\mid i{\in} I\}$.}
Let $S([a_i])$ be the symmetric group on the set $[a_i]=\{a\in
L\mid{\varepsilon}(a)=a_i\}$ and $S'([a_i])=\{u_i\in S([a_i])\mid
u_i(a_i)=a_i\}$ (of course, $S'([a_i])$ is a subgroup of
$S([a_i]))$, $i\in I.$ Denote by $\xmare{i\in I}{}S'([a_i])$ the
direct product of the groups $S'([a_i])$, $i\in I$. For every
$(u_i)_{i\in I}\in\xmare{i\in I}{}S'([a_i]),$ we construct an
\slat\ automorphism $f$ of $L$ by $f\mid \Fixe=1_{\Fixe}$ and
$f\mid[a_i]=u_i$, $i\in I$. In this way we defined a map
${\varphi}:\xmare{i\in I}{}S'([a_i])\longrightarrow {\rm
Aut}_{\cal E}(L).$ Moreover, note that ${\varphi}$ is a group
monomorphism.

On the other hand,  by the previous corollary, every element $f\in{\rm Aut}_{\cal E}(L)$ induces a lattice automorphism $f\mid\Fixe\in{\rm Aut}(\Fixe).$ Thus we defined another map $\psi:{\rm Aut}_{\cal E}(L)\longrightarrow {\rm Aut}(\Fixe)$, which is a group  epimorphism. Also, it is easy to see that $\Im\,{\varphi}=\Ker\,\psi$, therefore we have proved the next result.

\bigskip\n{\bf Proposition 2.} {\it With the above notations, there exists an exact sequence:
$$1\longrightarrow \xmare{i\in I}{}S'([a_i])\bldl{\varphi}{\rm Aut}_{\cal E}(L)\bldl\psi{\rm Aut}(\Fixe)\longrightarrow 1.$$}

Remark that we identified an important normal subgroup of ${\rm Aut}_{\cal E}(L)$:
$${\rm Aut}^0_{\cal E}(L)=\left\{f\in{\rm Aut}_{\cal E}(L)\mid  f|\Fixe=1_{\Fixe}\right\},$$
which is isomorphic to the direct product $\xmare{i\in I}{}S'([a_i]).$ A case when the group ${\rm Aut}_{\cal E}(L)$ itself is isomorphic to $\xmare{i\in I}{}S'([a_i])$ is described by the following corollary.

\bigskip\n{\bf Corollary 1.} {\it Under the same notations as in Proposition $2$, if the group ${\rm Aut}(\Fixe)$ is trivial, then we have:
$${\rm Aut}_{\cal E}(L)\cong\xmare{i\in I}{}S'([a_i]).$$}

There exist  many situations in which the group of all automorphisms of a lattice is trivial. One of them is obtained when the lattice is finite and fully ordered.

\bigskip\n{\bf Corollary 2.} {\it Let $\,(L,\wedge_{{\varepsilon}},\vee_{{\varepsilon}})\,$ be a canonical \slat\ having a finite fully ordered lattice of fixed points $\Fixe=\{a_1,a_2,...,a_n\}.$ Then the following group isomorphism holds:
$${\rm Aut}_{\cal E}(L)\cong\xmare{i=1}nS'([a_i]).$$}

Moreover, if $L$ itself is finite, we can estimate its number of \slat\ automorphisms.

\bigskip\n{\bf Corollary 3.} {\it Let $\,(L,\wedge_{{\varepsilon}},\vee_{{\varepsilon}})\,$ be a finite canonical \slat\ having a   fully ordered lattice of fixed points $\Fixe=\{a_1,a_2,...,a_n\}.$ If $m_i=|[a_i]|,$ $i=\overline{1,n}$, then the following equality holds:
$$\left|{\rm Aut}_{\cal E}(L)\right|=\prod_{i=1}^n(m_i-1)!.$$}

\subsection{Isomorphisms of subgroup \slats}

In this section we investigate isomorphisms between subgroup \slats.

Let $G_1,G_2$ be two groups and $L(G_1),L(G_2)$ be their subgroups
lattices. Remind that $G_1$ and $G_2$ are called
{\it$L$-isomorphic} if $L(G_1)\cong L(G_2)$ (a lattice isomorphism
from $L(G_1)$ to $L(G_2)$ will be called an {\it$L$-isomorphism}).
Now, let us consider the \slat\ structure on $L(G_i)$, $i=1,2.$ As
we have seen above, in this situation the lattices of fixed points
associated to $L(G_1)$, $L(G_2)$ are the normal subgroup lattices
$N(G_1),N(G_2)$ of $G_1$ and $G_2$, respectively. We say that
$G_1$ and $G_2$ are {\it${\cal E} L$-isomorphic} if the \slats\
$L(G_1)$ and $L(G_2)$ are isomorphic (an \slat\ isomorphism from
$L(G_1)$ to $L(G_2)$ will be called an {\it${\cal E}
L$-isomorphism}).

Our first goal is to establish some connections between these different types of isomorphisms from $G_1$ to $G_2$. Clearly, if $G_1\cong G_2$, then $G_1,G_2$ are both $L$-isomorphic and ${\cal E} L$-isomorphic. Also, if $G_1$ and $G_2$ are ${\cal E} L$-isomorphic, then, by Proposition 1, 2.1, $N(G_1)\cong N(G_2)$, but they are not necessarily $L$-isomorphic. Conversely, $N(G_1)\cong N(G_2)$ does not imply that $G_1$ and $G_2$ are ${\cal E} L$-isomorphic (for example, take $G_1$ the quaternion group and $G_2$ the dihedral group of order 8). Moreover, even the lattice isomorphism $L(G_1)\cong L(G_2)$ does not assure that $G_1$ and $G_2$ are ${\cal E} L$-isomorphic (for example, take $G_1$ a finite elementary abelian $p$-group and $G_2$ the nonabelian $P$-group which is $L$-isomorphic to $G_1$ (see [7], page 11)).\bigskip

The following result indicates us some classes of groups which are preserved by ${\cal E} L$-isomorphisms.

\bigskip\n{\bf Proposition 1.} {\it Let $G_1$ and $G_2$ be  two ${\cal E} L$-isomorphic groups. If $G_1$ is a simple group or a Dedekind $($in particular abelian$)$ group, then $G_2$ is also simple or Dedekind, respectively.}\bigskip

\n{\bf Proof.} Let $f:L(G_1)\longrightarrow  L(G_2)$ be an ${\cal
E} L$-isomorphism.

If $G_1$ is simple, then $|N(G_1)|=2.$ Since $N(G_1)\cong N(G_2)$, it results that $|N(G_2)|=2$ and hence $G_2$ is simple, too.

If $G_1$ is a Dedekind group, then any subgroup $H_1\in L(G_1)$ is normal in $G_1$ and its congruence class $[H_1]_1$ consists only of $H_1$. Because $f$ maps normal subgroups into normal subgroups and induces an one-to-one and onto map between the congruence classes $[H_1]_1$ and $[f(H_1)]_2$ of the \slats\ $L(G_1)$ and $L(G_2)$, respectively, it obtains $[H_2]_2=\{H_2\},$ for all $H_2\in N(G_2)$. This implies that every subgroup of $G_2$ is normal and hence $G_2$ is also Dedekind.\bigskip

Note that, for two groups of the above types, we are able to indicate some necessary and sufficient conditions in order to be ${\cal E} L$-isomorphic. In this way, two finite simple groups are ${\cal E} L$-isomorphic iff they have the same number of subgroups and two Dedekind (in particular abelian) groups are ${\cal E} L$-isomorphic iff they are $L$-isomorphic.\bigskip

Next, we shall present a property satisfied by ${\cal E} L$-isomorphisms between two finite  groups in the case when one of them is nilpotent.

\bigskip\n{\bf Proposition 2.} {\it Let $G_1,G_2$ be two finite groups, $\Phi(G_1),\Phi(G_2)$ be their Frattini subgroups and $f:L(G_1)\longrightarrow  L(G_2)$ be
an
 ${\cal E} L$-isomorphism. If $G_1$ is nilpotent, then $f(\Phi(G_1))\supseteq \Phi(G_2).$}\bigskip

 \n{\bf Proof.} Let $M_1$ be a maximal subgroup of $G_1$. Since $G_1$ is nilpotent, it follows that $M_1$ is normal in $G_1$ and so $M_2=f(M_1)$ is a normal subgroup of $G_2$. Let $H_2\in L(G_2)$ such that $M_2\subseteq H_2\subseteq G_2$ and assume that $H_2\ne G_2.$ Then $M_2$ is contained in the core $C_2$ of $H_2$ in $G_2$. As $f$ induces a lattice isomorphism between $N(G_1)$ and $N(G_2)$, it results that $C_1=f^{-1}(C_2)$ is a normal subgroup of $G_1$ and $M_1\subseteq C_1\subseteq G_1$. By the maximality of $M_1$, it obtains $M_1=C_1$ and therefore $M_2=C_2.$ This shows that $H_2\in[M_2]_2.$ But $f$ induces also an one-to-one and onto map between $[M_1]_1=\{M_1\}$ and $[M_2]_2$, thus $H_2=M_2.$ Because $f$ maps the maximal subgroups of $G_1$ into maximal subgroups of $G_2$, we have $f(\Phi(G_1))\supseteq\Phi(G_2).$

\bigskip\n{\bf Corollary.} {\it Let $G$ be a finite nilpotent group, $\Phi(G)$ be its Frattini subgroup and $f$ be an ${\cal E} L$-automorphism of $G$. Then $\Phi(G)$ is a fixed point of $f$.}\bigskip

\n{\bf Proof.} By Proposition 2, we have $f(\Phi(G))\supseteq \Phi(G)$. On the other hand,   applying Proposition 2 to $f^{-1}$, it results that $f^{-1}(\Phi(G)))\supseteq\Phi(G)$ and therefore $\Phi(G)\supseteq f(\Phi(G)).$ Hence the equality $f(\Phi(G))=\Phi(G)$ holds.

\bigskip\n{\bf Remark.} Assume that $G_1,G_2$ are two finite groups and let $f:L(G_1)\longrightarrow  L(G_2)$ be an ${\cal E} L$-isomorphism. By a well-known  result of H. Heineken (see [4]), under the additional conditions that $G_1$ is a noncyclic $p$-group and the derived subgroup of $G_2$ is nilpotent, it obtains that $G_2$ is also a $p$-group of the same order as $G_1$. In this case $f$ maps any principal series of $G_1$ into a principal series of $G_2$ and induces an one-to-one and onto map between the sets of maximal subgroups of $G_1$ and $G_2$. By Proposition 2, we have $f(\Phi(G_1))=\Phi(G_2)$ and thus $|\Phi(G_1)|=|\Phi(G_2)|.$ It follows that the vector spaces (over $F_p)$ $G_1/\Phi(G_1)$ and $G_2/\Phi(G_2)$ have the same dimension. Hence $G_1$ and $G_2$ can be generated by exactly the same number of generators.
\bigskip

Let $G_1,G_2$ be two groups, $f:L(G_1)\longrightarrow  L(G_2)$ be
an ${\cal E} L$-isomorphism and $H_1$ be a (normal) subgroup of
$G_1$. Since a normal subgroup of $H_1$ is not necessarily maped
by $f$ into a normal subgroup of $f(H_1)$, $f$ induces not an
${\cal E} L$-isomorphism between $H_1$ and $f(H_1).$ The situation
is different with respect to the factor groups of our two groups,
as shows the following lemma.

\bigskip\n{\bf Lemma.} {\it If  $G_1,G_2$ are two groups, $f:L(G_1)\longrightarrow  L(G_2)$ is an ${\cal E} L$-isomorphism and $H_1$ is a normal subgroup of $G_1$, then the map $$\bar f:L(G_1/H_1)\longrightarrow  L(G_2/f(H_1))$$ defined by $\bar f(K_1/H_1)=f(K_1)/f(H_1),$ for all $K_1/H_1\in L(G_1/H_1)$, is also an ${\cal E} L$-isomorphism.}\bigskip

In the hypothesis of the  above lemma consider $H_1$ be the
derived subgroup $D(G_1)$ of $G_1$. Then the groups $G_1/D(G_1)$
and $G_2/f(D(G_1))$ are ${\cal E} L$-isomorphic. But $G_1/D(G_1)$
is abelian and an ${\cal E} L$-isomorphism maps abelian groups
into Dedekind groups (see Proposition 1, 2.1), therefore
$G_2/f(D(G_1))$ is a Dedekind group $L$-isomorphic to
$G_1/D(G_1)$. It is well-known (for example, see [7], Theorem 6,
page 39) that a primary hamiltonian group cannot be $L$-isomorphic
to an abelian group. Thus, under a supplementary condition of type
$$\mbox{every  hamiltonian quotient of $G_2$ is primary,}\leqno(*)$$
it follows that $G_2/f(D(G_1))$ is also abelian and therefore
$D(G_2)\subseteq f(D(G_1))$ (mention that the author has not be
able to decide if without a condition of type $(*)$ it obtains the
commutativity of $G_2/f(D(G_1))).$ Hence we have proved the next
proposition.

\bigskip\n{\bf Proposition 3.} {\it Let $G_1,G_2$ be two groups, $D(G_1),D(G_2)$ be their derived subgroups and $f:L(G_1)\longrightarrow  L(G_2)$ be an ${\cal E} L$-isomorphism. If $G_2$ satisfies the condition $(*)$, then $f(D(G_1))\supseteq D(G_2).$}
\bigskip

By Proposition 3, we can easily see that the following result holds.

\bigskip\n{\bf Corollary.} {\it Let $G$ be a group which satisfies the condition $(*)$, $D(G)$ be its derived subgroup and $f$ be an ${\cal E} L$-automorphism of $G$. Then $D(G)$ is a fixed point of $f$.}\bigskip

As we have already seen in 2.1, an important normal subgroup of the group ${\rm Aut}_{\cal E}(L(G))$ associated to a group $G$ is ${\rm Aut}^0_{\cal E}(L(G))=\{f\in{\rm Aut}_{\cal E}(L(G))\mid f|N(G)=1_{N(G)}\}.$ We finish this section by indicating another two remarkable subgroups of ${\rm Aut}_{\cal E}(L(G)):$
\begin{itemize}
\item[--] the subgroup ${\rm Aut}^1_{\cal E}(L(G))$ consisting of all $L$-automorphisms of $G$ induced by group automorphisms;
\item[--] the subgroup ${\rm Aut}^2_{\cal E}(L(G))$ consisting of all $L$-automorphisms of $G$ of the type $f_a:L(G)\longrightarrow  L(G)$, $f_a(H)=H^a,$ for all $H\in L(G)$ $(a\in G)$.
\end{itemize}

Note that ${\rm Aut}^2_{\cal E}(L(G))$ is a normal subgroup of ${\rm Aut}^1_{\cal E}(L(G))$ and, also, it is contained in ${\rm Aut}^0_{\cal E}(L(G)).$

Finally, we find these subgroups in two situations.

\bigskip\n{\bf Examples.}
\begin{itemize}
\item[1)] For the symmetric group $S_3$ of degree 3, we have:\smallskip\\
${\rm Aut}^0_{\cal E}(L(S_3))={\rm Aut}^1_{\cal E}(L(S_3))={\rm Aut}^2_{\cal E}(L(S_3))\cong S_3.$
\item[2)] For the dihedral  group $D_4$ of order 8, we have:\smallskip\\
${\rm Aut}^0_{\cal E}(L(D_4))\cong S_4,\ {\rm Aut}^1_{\cal E}(L(D_4))\cong D_4,\ {\rm Aut}^2_{\cal E}(L(D_4))\cong \Z_2.$
\end{itemize}
\bigskip

\section*{References}
\begin{itemize}
\item[{[1]}] Birkhoff, G., {\it Lattice theory}, Amer. Math. Soc., Providence, R.I., 1967.
\item[{[2]}] Curzio, M.,  {\it Una caratterizzazione reticolare dei gruppi abeliani}, Rend. Mat. Appl. (7) 24 (1965), 1-10.
\item[{[3]}] Gr\"atzer, G., {\it General lattice theory}, Academic Press, New York, 1978.
\item[{[4]}] Heineken, H., {\it \"Uber die Charakterisierung von Gruppen durch gewisse Unter-gruppenverb\"ande}, J. Reine Angew. Math. 220 (1965), 30-36.
\item[{[5]}] Schmidt, R., {\it Subgroup lattices of groups}, de
Gruyter Expositions in Mathe\-ma\-tics 14, de Gruyter, Berlin, 1994.
\item[{[6]}] Suzuki, M.,  {\it Group theory,} I, II, Springer-Verlag, Berlin, 1982, 1986.
\item[{[7]}] Suzuki, M.,  {\it Structure of a group and
the structure of its lattice of subgroups}, Springer-Verlag, Berlin, 1956. \item[{[8]}] \c{S}tef\u anescu, M.,  {\it Introduction to group theory} (Romanian), Ed. Univ. "Al.I. Cuza" Ia\c si, 1993.
\item[{[9]}] T\u arn\u auceanu, M.,  {\it Actions of finite groups on lattices}, Seminar Series in Mathematics, Algebra 4, Univ. "Ovidius",
Constan\c ta, 2003.
\item[{[10]}] T\u arn\u auceanu, M.,  {\it ${\cal E}$-lattices}, accepted for publication in Italian Journal of Pure and Applied Mathematics.
\end{itemize}

\end{document}